\documentclass[11pt]{amsart}
\theoremstyle{plain}
\newtheorem{Thm}{Theorem}

\newtheorem{Pro}[Thm]{Proposition}

\begin{document}

%begin Topmatter
\title[liouville type theorems]
{Liouville type theorems for conformal Gaussian curvature
equation}

\author{Yihong Du, Li MA}

\address{Y.Du, Department of Mathematics \\
School of Sci. and TECHNOLOGY\\
University of New England \\
Armidale, NSW 2351\\
Australia} \email{ydu@turing.une.edu.au}

\address{Li Ma, Department of mathematical sciences \\
Tsinghua University \\
Beijing 100084 \\
China} \email{lma@math.tsinghua.edu.cn} \dedicatory{}
\date{August 6th, 2008}

\begin{abstract}

In this note, we study Liouville type theorem for conformal
Gaussian curvature equation (also called the mean field equation)
$$
-\Delta u=K(x)e^u, \; \; in \; R^2
$$
where $K(x)$ is a smooth function on $R^2$. When $K(x)=K(x_1)$ is a
sign-changing smooth function in the real line $R$, we have a
non-existence result for the finite total curvature solutions. When
$K$ is monotone non-decreasing along every ray starting at origin,
we can prove a non-existence result too. We use the methods of
moving planes and moving spheres.


\end{abstract}

\thanks{$^*$ The research is partially supported by the national research council of Australia,
the National Natural Science
Foundation of China 10631020 and SRFDP 20060003002 }
 \maketitle

\section{Introduction}

In this paper, we study the Liouville type results for conformal
Gaussian curvature equation (also called
 mean field equation or Lame-Emden equation)
\begin{equation}\label{mean}
-\Delta u=K(x)e^u, \; \; x=(x_1,x_2)\in R^2
\end{equation}
where is a smooth function on $R^2$. The geometrical meaning for
 the equation (\ref{mean}) is that the conformal metric $e^udx^2$ has its scalar curvature $K$.
 We firstly consider the case when $K(x)=K(x_1)$ is a sign-changing
smooth function in the real line $R$. This kind of problem arises
from the apriori bound for solutions via the blowing up argument.
In recent studies for the prescribed Gaussian curvature problem in
$R^2$ or the mean field equations, Radial symmetry, Liouville
theorems, and classification results for solutions with
\emph{finite energy} to equation (\ref{mean}) are obtained with
other assumptions on positivity or negativity of $K$, see
\cite{mc}, \cite{cni},\cite{ck},\cite{lin1}, \cite{li}, and
\cite{bt} therein for more references. In particular in \cite{mc},
the best existence result of solutions for a class of positive
functions $K$ has been obtained. In the work \cite{cni}, best
existence result for a class of negative functions $K$ has been
obtained. In the work \cite{lin1}, the behavior at infinity of
solutions with finite energy has been found when the function $K$
being positive or negative has a controlled polynomial growth.
When $K$ is non-positive, the quantity
$$
\alpha_1=\sup\{\alpha\in R;
\int_{R^2}|K(x)|(1+|x|^2)^{\alpha}dx<\infty\}
$$
plays an important role. When $K$ is positive with polynomial
growth, the total curvature
$$
\int_{R^2}K(x)e^{u(x)}dx
$$
plays the key role. One may see \cite{lin2} for more results.
However, there are relative few result for the case when the
function $K$ is changing sign or with no control of the growth. We
then show that Liouville type result is also true for a class of
positive radially monotone functions $K$. Since we have two kinds
of results with different assumptions of $K$, we state them
separately.

\emph{Results One}:

 We assume that $K(x)=x_1$ is non-trivial and $u$ is a smooth
 solution to (\ref{mean}) such that
$$
u(x)\leq u(0)=0, \; in \; \; R^2, \; \; \int_{R^2}|K(x)|
e^udx<+\infty.
$$
We remark that a smooth solution $u$ with $|K|e^u\in L^1(R^2)$ is
called \emph{finite energy} solution. We point out that one may
replace $K(x)=x_1$ by a nontrivial function with
$$
\partial_{x_1} K(x_1)\geq 0, \; in \; \; R^2
$$
in the result below.

One of the main purpose of this paper is to prove the following
Liouville type theorem.

\begin{Thm}\label{liouville}
Under the assumptions above and the finite total integral
$A:=\int_{R^2}x_1u^{u}>0$, there is no smooth solution to
(\ref{mean}).
\end{Thm}

Intuitionally, one may believe this true. The reason is that using
the Pohozaev identity (see (\ref{pohozaev}) in appendix) with the
vanishing boundary terms and $j=1$, we can immediately get
$$
0<\int_{R^2} e^u=0,
$$
which is absurd. Hence, Theorem \ref{liouville} is true by assuming
suitable decay conditions for boundary terms on large balls. Instead
of investigating this method, we shall use another method to prove
this result. Assume $u$ is a solution to (\ref{mean}). We shall
first derive an apriori estimate for solutions $u$ to (\ref{mean})
with finite total integral assumption. Then we use the moving plane
method to show that $u$ is monotone non-decreasing in $x_1$, which
leads a contradiction.

As an application of the Liouville theorem above, we may guess the
following
\begin{Thm}\label{thm1}
 Assume $\tilde{K}$ is a sign-changing smooth function on the two sphere $S^2$ such
that for small positive function $C_2>0$, there is a positive
constant $C_1>0$ satisfying $|\nabla \tilde{K}(y)|\geq C_1>0$ for
$$y\in\Gamma:=\{y\in S^2; |\tilde{K}(y)|\leq C_2\}.$$ Then there is
uniform constant $C>0$ such that
$$
|w(y)|\leq C, \; \; y\in \Gamma
$$
for all solutions $w$ to the conformal Gaussian curvature equation
$$
-\Delta_{S^2} w+2=\tilde{K}(y)e^w, \; \; y=(y_1,y_2,y_3)\in S^2
$$
\end{Thm}

This result has been proved by Chen and Li in \cite{Chenwl}. So we
shall not prove it, but we outline our formal argument of it. We
remark that the simpler a priori estimate for solutions on the
negative part of the function $K$ is given in \cite{ma08}. So we
consider only the result near zero set of the function $K$. Choose
$p\in S^2$ such that $K(p)<0$ and $K(-p)>0$. Let $x_1=y_1/(1-y_3),
x_2=y_2/(1-y_3)$ be the inverse stereographic projection mapping
$-p$ into origin in $R^2$. Then the spherical metric can be written
as
$$
\frac{4}{(1+|x|^2)^2}dx^2.
$$
Let $u(x)=w(y(x))+2\log\frac{2}{1+|x|^2}$ and let
$K(x)=\tilde{K}(y(x))$. Then we have
$$
-\Delta u=K(x)e^u, \; \; in \; R^2.
$$
Note that $$
\int_{S^2}\tilde{K}(y)e^{w(y)}d\sigma_y=\int_{R^2}\tilde{K}(y(x))e^{w(y(x))}\frac{4dx}{(1+|x|^2)^2}
=\int_{R^2}K(x)e^{u(x)}dx
$$
and
$$
\int_{S^2}|\tilde{K}(y)|e^{w(y)}d\sigma_y=\int_{R^2}|K(x)|e^{u(x)}dx.
$$
Hence, the conformal Gaussian curvature equation is reduced into
(\ref{mean}). Now, for any finite energy solution sequence $(u_j)$
to (\ref{mean}) with $\gamma_j=\{x; u_j(x)=0\}$, $u_j(x_j)\to
+\infty$ and $dist(x_j,\gamma_j)\to 0$, let
$d_j(x)=dist(x,\gamma_j)$.

 One can show that for some constant $\alpha$, $u_j(x)-\alpha \log
 d_j(x)$ is bounded in $\Gamma_j$ (see also \cite{Chenwl}).
Let
$$
v_j(x)=u_j(x+d_j(x)y)-\alpha \log d_j(x).
$$
Then we have a sub-convergence sequence, still denoted by $v_j$,
such that $v_j\to V$ in $C_{loc}^2(R^2)$ where $V$ satisfying
$$
-\Delta V=x_1e^V, \; \; x=(x_1,x_2)\in R^2.
$$
Using Theorem \ref{liouville} we find a contradiction. Hence, we
have proved Theorem \ref{thm1}.

We remark that similar result for scalar curvature problem has
been found by Chen and Li, C.Chen and C.S.Lin  (see \cite{CL97}
for more references). One may weaken the assumption on $K$ by
allowing $K$ to have large zero set, see \cite{DL}.

 \emph{Result Two}:

Using the moving sphere method we can prove
\begin{Thm}\label{thm2}
Let $K$ be a non-trivial positive $C^1$ function in $R^2$. Assume
that $K$ is non-decreasing along each ray $\{t\xi; t\geq 0\}$ for
every unit vector $\xi\in S^1$ with $x\cdot\nabla K(x)<2K(x)$ on
$R^2$ and
$$
K(\infty)=\lim_{|x|\to\infty}K(x)>0.
$$
Then the equation (\ref{mean}) has no smooth solution with
$\int_{R^2}Ke^udx<+\infty$.
\end{Thm}

We remark that if we have a solution in Theorem \ref{thm2}, then by
Cohn-Vossen inequality, we have \begin{equation}\label{vossen}
\int_{R^2} K(x)e^udx\leq 4\pi.
\end{equation}
We shall use this fact in our contrary argument.

The plan of the paper is below. We give some asymptotic behavior
estimate in section \ref{sect1}. In section \ref{sect2}, we prove
the Liouville theorem. In section \ref{sect3}, we prove Theorem
\ref{thm2}.

\section{Asymptotic behavior}\label{sect1}

To make the moving plane method get started at infinity, we need
to know the behavior of solutions at infinity. Let $u$ be a
solution to (\ref{mean}) with the upper bound $u(x)\leq u(0)$.

Let
$$
K(x)=K_+(x)-K_-(x),
$$
where $K_+$ be the positive part of $K$ and $K_-$ is the negative
part of $K$.

Define $v=v_0-v_1$, where
$$
v_0(x)=\frac{1}{2\pi}\int_{R^2}\log \frac{|x-y|}{|y|}
K_+(y)e^{u(y)}dy
$$
and
$$
v_1(x)=\frac{1}{2\pi}\int_{R^2}\log \frac{|x-y|}{|y|}
K_-(y)e^{u(y)}dy
$$

 Let $g_0(x)=K_+(x)e^{u(x)}$. Let $R>0$ be a large number. We write
$$
R^2_+=\{x\in R^2; x_1\geq 0\}=B_R^+\bigcup T_1\bigcup T_2,
$$
where $B_R^+(0)=\{x\in R^2_+; |x|\leq R\}$,
$$ T_1=\{y=(y_1,y_2); y_1>0, |y-x|\leq \frac{|x|}{2}\},
$$
and
$$
T_2=\{y=(y_1,y_2); y_1>0, |y-x|\geq \frac{|x|}{2}\; and \; |y|\geq
R
$$

Then for $x_1>> 1$, we have
$$
2\pi
v_0(x)=(\int_{B_R^+}+\int_{T_1}+\int_{T_2})\log\frac{|x-y|}{|y|}g_0(y)dy
\equiv I_0+I_1+I_2.
$$

Note that for $y\in T_1$, we have $|y|\leq |x|/2$ and $|x-y|\leq
|y|$. Hence, $I_1\leq 0$.

It is clear that
$$
I_0\leq \log|x|\int_{B_R^+} g_0(y)dy+C
$$

For $y\in T_2$, we have $|y|\geq R >1$ and $|x-y|\leq |x|+|y|\leq
|x||y|$. Then we have
$$
I_2\leq \log|x|\int_{T_2} g_0(y)dy.
$$

Hence, we have
$$
2\pi v_0(x)\leq I_0+I_2\leq \log |x|(\int_{B_R^+}+\int_{T_2})
K_+(y)e^{u(y)}dy +C.
$$

It is also clear that for $x_1<<-1$
$$
2\pi v_0(x)=(\int_{B_R^+}K_+(y)e^{u(y)}dy+\circ(1))\log |x|.
$$
The lower bound for $v_0$ can be obtained in the same way (see
also the treatment about $v_1$ below).

We now find the the lower bound for $v_1$. Define
$$
R^2_-=\{x\in R^2; x_1\leq 0\}=B_R^+\bigcup S_1\bigcup S_2,
$$
where $B_R^-(0)=\{x\in R^2_-; |x|\leq R\}$,
$$ S_1=\{y=(y_1,y_2); y_1<0, |y-x|\leq \frac{|x|}{2}\},
$$
and
$$
S_2=\{y=(y_1,y_2); y_1<0, |y-x|\geq \frac{|x|}{2}\; and \; |y|\geq
R
$$

Then we have that for $x_1<<-1$,
$$
2\pi v_1(x)\geq \log |x|(\int_{B_R^-}+\int_{S_2}) K_-(y)e^{u(y)}dy
+C
$$
and $x_1>>1$
$$
2\pi v_1(x)=(\int_{B_R^-}K_-(y)e^{u(y)}dy+\circ(1))\log |x|+C.
$$

In the same way, we have the same control in $x_2$ direction.

As for the upper bound of $v_1(x)$ for $|x|$ large, noticing
$u(x)\leq 0$, we may invoke the Lemma 2.2 and Lemma 2.3 in
\cite{lin1}. For completeness, let us do it here.

Let $g_1(x)=K_-(x)e^{u(x)}$. As before, for $|x|>> 1$, we write
$$
2\pi
v_1(x)=(\int_{B_R^+}+\int_{S_1}+\int_{S_2})\log\frac{|x-y|}{|y|}g_1(y)dy
\equiv I_0+I_1+I_2.
$$
For $y\in S_2$, we have $|x-y|\leq |x|+|y|\leq |x||y|$ and
$$
I_2\leq \log \int_{S_2}g_1(y)dy.
$$

It is also easy to see that
$$
\log |x|\int_{B_R^+}g_1(y)dy-C\leq I_0\leq \log
|x|\int_{B_R^+}g_1(y)dy+C.
$$

Note that $|y|\geq |x|/2\geq|x-y|$ in $S_1$, and we have
$$
I_1= (\int_{|x-y|\leq |x|^{-4/3}}+\int_{|x|^{-4/3}\leq |x-y|\leq
|x|/2})\log\frac{|x-y|}{|y|}g_1(y)dy.
$$

By computation, we have
$$
\int_{|x-y|\leq |x|^{-4/3}}|\log\frac{1}{|x-y|}|\leq |x|^{-1/2},
$$
and
$$
|\int_{|x|^{-4/3}\leq |x-y|\leq
|x|/2}\log\frac{|x-y|}{|y|}g_1(y)dy|\leq 8\log |x|\int_{|x-y|\leq
|x|/2} g_1(y)dy.
$$
Using $u(x)\leq 0$, we have
$$
|I_1|\leq \epsilon \log|x|+C
$$
for $|x|$ large. Hence, we have
$$
2\pi v_1(x)\leq (\int_{B_R+T_2} g_1(y)dy+\epsilon)\log|x|+C.
$$

We remark that for the lower bound of $v_1$, we use
$$
I_0\geq \log |x|\int_{B_R^+}g_1(y)dy-C
$$
and for large $R>>1$,
$$
I_2\geq -1
$$
where we have used the fact $|y|\leq 4|x-y|$, which implies that
$$0\leq \log 4+\log \frac{|x-y|}{|y|}$$ for $y\in S_2$ and $$
I_2\geq -\log 4\int_{S_2} g_1(g)dy.
$$

Putting all these together, we obtain

\begin{Pro}\label{du1} Assume that
$$
A:=\int_{R^2}x_1e^{u}>0.
$$
Then there is a positive constant $A$ such that
 for any $\epsilon>0$, there exists a large constant $R(\epsilon)>0$, for $|x|\geq
 R(\epsilon)$, it holds
 $$
[A-\epsilon]\log|x| -C\leq 2\pi v(x)\leq [A+\epsilon]\log|x| +C.
 $$
\end{Pro}

Recall the following well-known Liouville type theorem for harmonic
functions.
\begin{Thm}\label{chenlin}
Assume that $w$ is a harmonic function in $R^2$ such that
$$
w(x)\leq A\log |x|+C
$$
for $|x|\geq R_0$, where $R_0$ is some fixed positive function. Then
the harmonic function  $w(x)$ is a constant.
\end{Thm}

By construction, we have $$ \Delta (u+v)=0, \; \; x\in R^2.
$$
Recall that $u$ is bounded from above, so we have
$$
u(x)+v(x)\leq A_0\log |x|+C
$$
at infinity for some $A_0\geq 0$. Then we can use the Liouville
Theorem above for harmonic functions to show that $u+v$ is constant.

\begin{Thm}\label{apriori}
$$u(x)+v(x)=Constant.
$$
\end{Thm}

\emph{Remark:} Using Theorem \ref{apriori}, we also ave that
$$
u(x)\geq -\alpha \log|x|+C.
$$

\section{moving plane method }\label{sect2}

We shall use the moving plane method to prove  that
$\partial_{x_1}u>0$. Using the fact that $u(0)=0$, we then get
$u(x_1,0)=0$ for all $x_1\geq 0$, which is a contradiction to the
property that $\lim_{x_1\to \infty} u(x_1,0)=-\infty$. Then we have
proved Theorem \ref{liouville}.

In doing the moving plane method, we let for any real $\lambda$,
$$
T_{\lambda}=\{x; x_1=\lambda\}, \; \; \Sigma_{\lambda}=\{x;
x_1<\lambda\}
$$
and
$$
x^{\lambda}=(2\lambda-x_1,x_2).
$$

Let
$$
w_{\lambda}=u_{\lambda}(x)-u(x)
$$
where
$$
u_{\lambda}(x)=u(x^{\lambda}):=u(2\lambda-x_1,x_2).
$$
Then we have
$$
\Delta w_{\lambda}(x)=K(x)e^u-K(x^{\lambda})e^{u_{\lambda}}.
$$

\emph{Claim: }
\begin{equation}\label{plane}
w_{\lambda}(x)>0, \; \; for \; x\in \Sigma_{\lambda}\; \; and \;
\lambda\in \mathbb{R}.
\end{equation}

We shall prove this Claim in two steps.

\emph{Step one. }(\ref{plane}) is true for $x\in \Sigma_{\lambda}$
and $\lambda\leq 0$.

In this step, we first show the following.

\emph{Assertion:}
$$
\Delta w_{\lambda}(x)<0, \; whenever \; \; w_{\lambda}(x)\leq 0 \;
for \; x\in \Sigma_{\lambda}, \lambda\leq 0.
$$

Suppose that $w_{\lambda}(x)=u_{\lambda}(x)-u(x)\leq 0$ for some
$x\in \Sigma_{\lambda}$. If $x_1^{\lambda}\geq 0$, then
$K(x^{\lambda})\geq 0$. Hence,
$$
\Delta w_{\lambda}(x)=K(x)e^u-K(x^{\lambda})e^{u_{\lambda}}<0.
$$
If $x_1^{\lambda}< 0$, then $|x^{\lambda}_1|\leq |x_1|$ and
$$
\Delta
w_{\lambda}(x)=-|K(x)|e^u+|K(x^{\lambda})|e^{u_{\lambda}}\leq0.
$$
Note that near infinity $|x|=+\infty$, by using Theorem
\ref{apriori} and Proposition \ref{du1}, $w(x)$ is bounded by
$2\epsilon\log |x|$.

 We now let
$$
g_{\lambda}(x)=\log(\lambda+3-x_1)+\log ((\lambda+3-x_1)^2+x_2^2),
\; \; x\in \Sigma_{\lambda}.
$$

Note that
$$
\Delta g_{\lambda}(x)=-(\lambda+3-x_1)^{-2}<0
$$
and near $|x|=\infty$, $g_{\lambda}(x)$ has the behavior no less
than $2\log |x|$.

Let
$$
\bar{w}(x)=\frac{w_{\lambda}(x)}{g_{\lambda}(x)},
$$
which is bounded near $|x|=\infty$. Let $x_0$ be the point such that
$$\bar{w}(x_0):=\inf\bar{w}(x).$$ Such an $x_0$ can be found since
$|\bar{w}(x)|$ is arbitrary small (since $\epsilon$ can be small).
 Consider the function
$$
\Delta\bar{w}(x)
$$
at $x_0$. Using $$ \nabla \bar{w}(x_0)=0, \; \;
\Delta\bar{w}(x_0)\geq 0,
$$
we obtain that
$$
0\geq \Delta w_{\lambda}(x_0)=g_{\lambda}(x_0)\Delta
\bar{w}(x_0)+\bar{w}(x_0)\Delta g_{\lambda}(x_0)
>0,
$$
Which is absurd.

Let
$$
\lambda_0=\sup\{\lambda; w_{\mu}(x)>0, x\in\Sigma_{\mu}\; for \;
\mu<\lambda\}.
$$
Clearly $\lambda_0> 0$.

\emph{Step two. } $\lambda_0=+\infty$.

For otherwise, we assume
\begin{equation}\label{absurd}
\lambda_0<+\infty. \end{equation}
 By definition, we have a sequence
$\lambda_j>\lambda_0$ with $\lim \lambda_j\to \lambda_0$ and
$$
\inf_{\Sigma_{\lambda_j}}w_{\lambda_j}(x)<0.
$$
Write
$$
w_j(x)=w_{\lambda_j}(x).
$$
As in the Step one, we want to show that for large $j$,
\begin{equation}\label{keystep}
\Delta w_j<0, \;  \; whenever \; \; w_j(x)\leq 0, \; for \; x\in
\Sigma_{\lambda_{j}}.
\end{equation}
Once this is done, we can repeat the Step one to get a
contradiction to (\ref{absurd}).

We argue by contradiction again. Assume $x_j=(x_{j1},
x_{j2})\in\Sigma_{\lambda_{j}}$ such that $w_j(x_j)\leq 0$ and
$\Delta w_j(x_j)\geq 0$. for short we write by $\lambda=\lambda_j$.
Since $\lambda_j\geq \lambda_0\geq 0$, $K(x_j^{\lambda})>0$ and
$$
0\leq \Delta
w_j(x_j)=K(x_j)e^{u(x_j)}-K(x_j^{\lambda})e^{u_{\lambda_j}(x_j)}.
$$
Then we have $K(x_j)\geq 0$ and  $0\leq x_{j1}\leq \lambda_j$.

Recall that by continuity, we have
$$
w_{\lambda}(x)\geq 0, \; x\in \Sigma_{\lambda_0} $$ and
$$
\partial_{x_1}w_{\lambda}(x)\geq 0, \; for \; x_1=\lambda_0.
$$
Using the strong maximum principle and Hopf's boundary point
lemma, we have
$$
w_{\lambda}(x)> 0, \; x\in \Sigma_{\lambda_0} $$ and $$
\partial_{x_1}w_{\lambda}(x)>0, \; for \; x_1=\lambda_0.
$$

If $(x_j)$ is bounded, we may further assume that
$x_0=\lim_{j\to\infty} x_j$. Then either $x_0\in
\Sigma_{\lambda_0}$, which gives $w_{\lambda_0}(x_0)\leq 0$, or
$x_0\in T_{\lambda_0}$ which implies that $|\nabla
w_{\lambda_0}(x_0)|=0$. All these give a contradiction.

Hence  $(x_j)$ is unbounded, which implies that $x_{j2}$ is
unbounded. Let
$$
\phi_j(x)=u_j(x+(0,x_{j2}))-M_j
$$
where $M_j=u_j(0,x_{j2})\to -\infty$. Then we have that $\phi_j$
is locally bounded from above satisfying
$$
-\Delta \phi_j= x_1 e^{M_j}e^{\phi_j}, \; in \; R^2.
$$

Using the Harnack inequality, we have a locally uniformly convergent
subsequence, still denoted by $\phi_j$ with  its limit $\phi$, which
satisfies
$$
-\Delta \phi=0, \; in \; R^2.
$$
Hence by Liouville theorem, $\phi=0$.

Since $\phi=0$, by locally uniformly convergence of $(\phi_j)$, we
have that for any $\epsilon>0$, there is a $j_0$ such that for
$j>j_0$,
$$
\max_{B((x_{j1},0),1+\lambda_0)}|\nabla u_j(x)|<\epsilon.
$$
Let $$ H(x)=x_1e^u-(2\lambda-x_1)e^{u_{\lambda}}.
$$ Then for any $x=(t, x_2)$ where $t\in [x_{j1},\lambda_j]$,
$$
\partial_{x_1}H(x)=e^u(1+x_1\partial_1u)
+e^{u_{\lambda}}(1+x^{\lambda}_1\partial_1u^{\lambda})>0.
$$
This gives us that for $x_1< \lambda$,
$$ H(x)<
H(\lambda,x_2)=0.
$$
 Hence,
$$
\Delta
w_j(x_j)=H(x_j)=(x_1e^u-(2\lambda-x_1)e^{u_{\lambda}})|_{x_j}<0,
$$
which yields a contradiction. Thus, (\ref{keystep}) is true. By
this we have $\lambda_0=+\infty$ and the Claim is true with
$\partial_{x_1}u>0$. In particular, $\partial_{x_1}u(0,0)>0$.

Since $u(0,0)=0$ is the maximum of $u$, we have
$$
\partial_{x_1}u(0,0)=0.
$$
A contradiction. Hence no such solution $u$ exists.

\section{moving sphere method }\label{sect3}

In this section, we prove Theorem \ref{thm2}.

\begin{proof}
Without loss of generality, we may assume that $K(\infty)=1$. We
shall use the method of moving spheres, which is a little bit easier
than the moving plane method since we only use the maximum principle
in bounded balls.

 Let
$$v(x)=u(\frac{x}{|x|^2})-4\log|x|.$$

For $\lambda>0$, $$
v_{\lambda}(x)=u({\lambda}^2\frac{x}{|x|^2})+4\log\lambda-4\log|x|.
$$
Then  $v_{\lambda}(x)$ satisfies
$$
\Delta v_{\lambda}+K(\lambda x/|x|^2)exp(v_{\lambda})=0.
$$

Let $w_{\lambda}(x)=v_{\lambda}(x) -u(x)$. Then $w_{\lambda}$
satisfies
$$
\Delta w_{\lambda}(x)+ b_{\lambda}(x)w_{\lambda}(x)=Q_{\lambda}(x)
$$
where
$$
b_{\lambda}(x)w_{\lambda}(x)=K(\lambda
x/|x|^2)(\frac{exp(v_{\lambda})-exp(u)}{v_{\lambda}-u})
$$
and
$$
Q_{\lambda}(x)=(K(\lambda x)-K(\lambda x/|x|^2))exp(u)
$$

By our monotone assumption on $K$, we have $Q_{\lambda}(x)\leq 0$
for $|x|<1$ and $\lambda>0$.

We claim that
\begin{equation}\label{monotone}
w_{\lambda}(x)>0, \; |x|<1\; \; and \; \lambda>0.
\end{equation}

Assume that $w_{\lambda}(x_0)=\inf_{B_1} w_{\lambda}<0$ for some
$x_0\in B_1$ and some small $\lambda\leq 1$, we have
$$
\Delta w_{\lambda}(x_0)\leq 0.
$$
But this is impossible by strong maximum principle. Hence
(\ref{monotone}) is true for small $\lambda>0$. Let
$$
\lambda_1=\sup \{\lambda; w_{\mu}(x)>0 \; for \; all \; 0<\mu\leq
\lambda\; and \; |x|<1\}.
$$

Then we have
$$
\lambda_1=+\infty.
$$
For otherwise, we have $\lambda_1<\infty$. Note that
$w_{\lambda_1}$ may have singularity in $x=0$. Using the maximum
principle to $w_{\lambda_1}\geq 0$ in $0<|x<1$, we have
$w_{\lambda_1}(x)>0$. Hence for any small $\epsilon>0$, there is
exists a $\delta>0$ such that for $|\lambda-\lambda_1|<\delta$, we
have
$$
\inf_{|x|=\epsilon}w_{\lambda}(x)>0.
$$
Arguing as above (as for small $\lambda$ before the definition of
$\lambda_1$), we have that for $|\lambda-\lambda_1|<\delta$ and
$0<|x|\leq \epsilon$,
$$
w_{\lambda}(x)>0, \; for \; 0<|x|\leq \epsilon.
$$

 By the definition of $\lambda_1$, we
have $\lambda_n>\lambda_1$ and $0\not=x_n\in B_1$ such that
$\lambda_n\to\lambda_1$ and
$$
w_{\lambda_n}(x_n)=\inf_{B_1}w_{\lambda_n}<0.
$$

With loss of generality, we can let $x_0=\lim_n x_n$. Then $|x_0|=1$
and $\nabla w_{\lambda_1}(x_0)=0$, which is a contradiction to
Hopf's boundary point lemma. Hence the Claim (\ref{monotone}) is
true. The Claim (\ref{monotone}) and the Hopf boundary lemma imply
that, for $y=x/|x|^2$,
$$
v(x/|x|^2)-2\log |x|=v(y)+2\log |y|
$$
is monotone decreasing along each ray, i.e., for $t>0$,
\begin{equation}\label{monotone2}
\partial_t(v(ty)+2\log t|y|)\leq 0.
\end{equation}

We now claim that there exists a uniform constant $C>0$ such that
\begin{equation}\label{upper}
u(x)+4\log |x|\leq C
\end{equation}
 in
$R^2$. For otherwise, we have a sequence $(x_j)\subset R^2$ such
that
$$
|x_j|\to \infty
$$
and
$$
u(x_j)+4\log |x_j|\to +\infty.
$$
That is saying that for $y_j=x_j/|x_j|^2$, $|y_j|\to 0$,
$$
v(y_j)\to +\infty.
$$
For $R>0$, define $$ S_j(y)=(R^2-|y-y_j|^2)e^{v(y)},
$$
in $B(y_j,R)$. Note that
$$
S_j(\bar{y}_j)=\max_{B(y_j,R)}S_j(y)\geq S_j(y_j)=exp(v(y_j)+2\log
R)\to +\infty.
$$
Define $M_j=-2\log\delta_j=v(\bar{y}_j)$ and
$$
\bar{v}_j(\eta)=v(\bar{y}_j+\delta_j y)-M_j.
$$
 Using the standard blowing up method (see \cite{li}) we
know that the renormalization $\bar{v}_j$  uniformly converges to
the bubble solution $U_0$ in $C_{loc}^2(R^2)$. Using the symmetry
of $U_0$ at $\eta_0$, we know that $\bar{v}_j$ has a local maximum
at $\eta_j$ with $\lim_j\eta_j=\eta_0$. Here $\eta_j$ is the point
such that $\bar{y}_j=y_j+\delta_j\eta_j$.  Going back to $v$, $v$
has a local maximum at $y*_j=y_j+\delta_j \eta_j$. Then we have at
$t=1$,
$$
\partial_t(v(ty*_j)+2\log t)\leq 0, \; and \; \; \partial_tv(ty*_j)=0,
$$
which are contrary each other.

Since $K({\infty})=1$, we may show that
\begin{equation}\label{lower} u(x)\geq -4\log |x|+c
\end{equation}
for large $|x|>>1$. In fact, there is another way. According to
Lemma 2.2 in \cite{lin1}, we have $$ u(x)\geq -{\beta}\log |x|+C_1
$$
for $\beta=\frac{1}{\pi}\int Ke^udx$ and $|x|$ large. By Cohn-Vossen
inequality (\ref{vossen}) we have $\beta\leq 4$. Then we have
$$ -4\log |x|+c\geq u(x)\geq -\beta\log |x|+C_1,
$$
which implies that $$ (4-\beta)\log |x|\leq C, \quad for\quad |x|>>1
$$
and then $\beta\geq 4$. Hence, $\beta=4$. Combining (\ref{upper})
with (\ref{lower}), we get
$$
u(x)=-4\log |x|+q(x)
$$
with $q(x)$ being bounded at $|x|\to \infty$. Applying the Pohozaev
identity in the appendix (see also \cite{mawei})
$$
-\int_{\partial B_R}[R(|\nabla
u|^2/2-|\partial_ru|^2-Ke^u)]=\int_{B_R}(K+\frac{1}{2}x\cdot \nabla
K)e^u
$$
to $u$ in a large ball $B_R$ and arguing as in pages 136-137 in
\cite{lin1} we then get
$$ \lim_{R\to \infty} \int_{B_R}x\cdot \nabla Ke^u=16\pi.
$$
By our assumption $x\cdot\nabla K<2K$, we get a contradiction. Then
we are done.
\end{proof}

\section{appendix}
In this appendix, we present the Pohozaev identity for equation
(\ref{mean}). This fact is well-known to experts (see
\cite{mawei}). Since the proof is different from the higher
dimensional case, we give a proof below.

Let $F(x,u)=K(x)e^u$. Integrating both sides of (\ref{mean}) over
$B_R(0)$, we have
$$
-\int_{\partial B_R(0)} \mu_ju_j=\int K e^u.
$$
Here $\int=\int_{B_R(0)}$ and $\mu_j=x_j/R$.

Multiplying both sides of (\ref{mean}) by
$r\partial_ru=x_j\partial_ju$ and integrating over $B_R(0)$, we
have
$$
\int (x_j\partial_ju)_iu_i-\int_{\partial
B_R(0)}x_ju_j\nu_iu_i=\int Kx_j\partial_je^u.
$$
 Note that
$$
\int (x_j\partial_ju)_iu_i=\int (|\nabla u|^2+x_ju_{ij}u_i)
$$
and
$$
\int x_ju_{ij}u_i=\frac{1}{2}\int x_j(|\nabla u|^2)_j=-\int
|\nabla u|^2+\frac{1}{2}\int_{\partial B_R(0)}x_j\nu_j|\nabla
u|^2.
$$
Then we have
$$
\int (x_j\partial_ju)_iu_i=\frac{1}{2}\int_{\partial
B_R(0)}x_j\nu_j|\nabla u|^2.
$$

Since
$$
\int Kx_j\partial_je^u=-\int (2Ke^u+x_j\partial_j
Ke^u)+\int_{\partial B_R(0)} x_j\nu_j Ke^u,
$$
we have
$$
\frac{1}{2}\int_{\partial B_R(0)}x_j\nu_j|\nabla
u|^2-\int_{\partial B_R(0)}x_ju_j\nu_iu_i-\int_{\partial B_R(0)}
x_j\nu_j Ke^u=-\int (2Ke^u+x_j\partial_j Ke^u).
$$
This is the standard \emph{Pohozaev identity} for the mean field
equation in $B_R(0)$. Sometimes, people would like to use another
form of it, which is
$$
\frac{1}{2}\int_{\partial B_R(0)}x_j\nu_j|\nabla
u|^2-\int_{\partial B_R(0)}x_ju_j\nu_iu_i-\int_{\partial B_R(0)}
x_j\nu_j Ke^u $$
$$-2\int_{\partial B_R(0)} \mu_ju_j=-\int
x_j\partial_j Ke^u.
$$

Similarly, by multiplying by $u_j$, we can get another Pohozaev
identity:
\begin{equation}\label{pohozaev}
\int\partial_jK e^u=\int_{\partial B_R(0)}
(\nu_ku_ku_j-\frac{1}{2}|\nabla u|^2\nu_j+Ke^u\nu_j)
\end{equation}

All these Pohozaev identities are useful in the study of mean field
equation.

\end{document}